\documentclass{amsart}
\usepackage[T1]{fontenc}
\usepackage[utf8]{inputenc}
\usepackage[english]{babel}
\usepackage{csquotes}

\usepackage{biblatex}
\renewbibmacro{in:}{%
        \ifentrytype{article}{}{\printtext{\bibstring{in}\intitlepunct}}}

\usepackage[hang,flushmargin]{footmisc} 

\usepackage{amsmath}
\usepackage{amssymb}
\usepackage{amsthm}

\usepackage{abstract}
\addto\captionsenglish{}

\newtheorem{defn}{Definition}[section]
\newtheorem{theo}[defn]{Theorem}
\newtheorem{prop}[defn]{Proposition}
\newtheorem{lemma}[defn]{Lemma}
\newtheorem{cor}[defn]{Corollary}

\newcommand{\N}{\mathbb{N}}

\newcommand{\func}[3]{#1:#2\rightarrow#3}
\newcommand{\norm}[2]{{||#1||}_{#2}}
\newcommand{\BL}{\textup{BL}}

\renewcommand{\epsilon}{\varepsilon}
\renewcommand{\theta}{\vartheta}
\renewcommand{\phi}{\varphi}

\title{Haar Null Closed and Convex Sets in Separable Banach Spaces}
\author{Davide Ravasini}
\address{Institut für Mathematik, Universität Innsbruck \newline\indent Technikerstraße 13, 6020 Innsbruck, Austria}
\email{davide.ravasini@uibk.ac.at}

\begin{document}
\maketitle
\let\thefootnote\relax\footnote{\today \newline \indent\emph{2020 Mathematics Subject Classification:} 46B20 (primary), 46B42, 52A07 (secondary). \newline
}

\begin{abstract}
\noindent \textsc{Abstract}. Haar null sets were introduced by J.P.R.\ Christensen in 1972 to extend the notion of sets with zero Haar measure to nonlocally compact Polish groups. In 2013, U.B.\ Darji defined a categorical version of Haar null sets, namely Haar meagre sets. The present paper aims to show that, whenever $C$ is a closed, convex subset of a separable Banach space, $C$ is Haar null if and only if $C$ is Haar meagre. We then use this fact to improve a theorem of E.\ Matou\v{s}kov\'{a} and to solve a conjecture proposed by Esterle, Matheron and Moreau. Finally, we apply the main theorem to find a characterisation of separable Banach lattices whose positive cone is not Haar null.
\end{abstract}

\section{Introduction}
\label{sec:intro}
Let $G$ be a Polish group. A Borel set $E\subset G$ is \emph{Haar null} if one can find a Borel probability measure $\mu$ on $G$ such that $\mu(gEh)=0$ for all $g,h\in G$. The measure $\mu$ is said to \emph{witness} that $E$ is Haar null and can always be assumed to have compact support, as every Borel probability measure on a Polish space is inner regular. This definition was given for the first time in Abelian Polish groups by J.P.R.\ Christensen in \cite{chris_1972} as a way to generalise the notion of sets with zero Haar measure in nonlocally compact Polish groups, where the Haar measure is not available. Indeed, Haar null sets and sets with zero Haar measure agree on Abelian, locally compact Polish groups. Moreover, the family of Haar null sets maintains the property of being a translation-invariant $\sigma$-ideal of the Borel $\sigma$-algebra of $G$. We say that a Borel set $E$ is \emph{Haar positive} or \emph{nonnegligible} if it is not Haar null. Drawing inspiration from Christensen's idea, U.B.\ Darji introduced in \cite{dar_2013} a categorical variant of the same concept. A Borel set $E$ in a Polish group $G$ is said to be \emph{Haar meagre} if there are a compact metric space $M$ and a continuous function $\func{f}{M}{G}$ such that $f^{-1}(gEh)$ is meagre in $M$ for every $g,h\in G$. The function $f$ is said to \emph{witness} that $E$ is Haar meagre. As in the case of Haar null sets, meagre sets and Haar meagre sets agree on locally compact Polish groups and the family of Haar meagre sets is a translation-invariant $\sigma$-ideal of the Borel $\sigma$-algebra of $G$. The reader is invited to have a look at the survey articles \cite{bog_2018}, \cite{en_2020} and at \cite{benlin}, Chapter 6, for more detailed information and results about Haar null and Haar meagre sets.

There is a sufficient, but rather demanding condition which ensures that a given Borel set $E$ in a Polish group $G$ is Haar positive and which can be expressed as follows: for every compact set $K\subset G$, there are an open set $V\subseteq G$ and $g,h\in G$ such that $K\cap V\ne\varnothing$ and $K\cap V\subset gEh$. Indeed, suppose that $E$ enjoys this property, let $\mu$ be any Borel probability measure on $G$ with compact support $K$ and find an open set $V\subseteq G$ and $g,h\in G$ such that $K\cap V\ne\varnothing$ and $K\cap V\subset gEh$. Clearly, $\mu(K\cap V)>0$, because $K$ is the support of $\mu$, hence $\mu(gEh)>0$. Sets with this property are called \emph{compactivorous}, can be defined in any Hausdorff topological group and were mentioned for the first time in \cite{emm_2016} by Esterle, Matheron and Moreau, who introduced them in separable Banach spaces. Compactivorous sets cannot be Haar meagre. Indeed, let $E$ be a Borel compactivorous set in a Polish group $G$, let $M$ be a compact metric space and let $\func{f}{M}{G}$ be a continuous function. Set $K=f(M)$ and find $g,h\in E$ such that $K\cap gEh$ has nonempty relative interior in $K$. Let $V\subset K\cap gEh$ be an open subset of $K$. Then $f^{-1}(V)$ is an open subset of $M$ such that $f^{-1}(V)\subset f^{-1}(gEh)$, hence $f^{-1}(gEh)$ is not meagre in $M$. Since $M$ and $f$ are arbitrary, it follows that $E$ is not Haar meagre. In the case of $\textup{F}_\sigma$ sets, the converse is true as well.

\begin{prop}
\label{prop:intro}
Let $E$ be an $\textup{F}_\sigma$ set in a Polish group $G$. Then $E$ is compactivorous if and only if it is not Haar meagre.
\end{prop}
\begin{proof}
We have already observed that if $E$ is compactivorous, then $E$ is not Haar meagre. Conversely, let ${\{C_n\}}_{n=1}^\infty$ be a countable family of closed sets such that 
\[ E=\bigcup_{n=1}^\infty C_n. \]
If we assume that $E$ is not compactivorous, then there is a compact set $K$ such that $K\cap gEh$ has an empty relative interior in $K$ for all $g,h\in G$. Let $\func{i}{K}{G}$ be the inclusion map. Observe that $K\cap g\,C_nh$ is a closed subset of $K$ with empty interior for every $g,h\in G$ and every $n$, which implies that it is nowhere dense. This shows that 
\[ K\cap gEh=\bigcup_{n=1}^\infty K\cap g\,C_nh \] 
is meagre for every $g,h\in G$, hence $i$ witnesses that $E$ is Haar meagre. 
\end{proof}

In particular, a closed set which is not Haar meagre must be Haar positive. What we want to prove is that the opposite implication holds for closed convex subsets of separable Banach spaces: given a closed, convex set $C$ in a separable Banach space, $C$ is Haar null if and only if $C$ is Haar meagre. 

In Section \ref{sec:C-fit}, the space of $C$-fit sequences is introduced. This will prove to be a fundamental tool in the development of the main argument, which is carried out in Section \ref{sec:hpccs}. Section \ref{sec:lattices} will present an application of the main theorem to the characterisation of separable Banach lattices with a nonnegligible positive cone. The standard notation of Banach space theory is used: if $X$ is a Banach space, we denote by $X^*$ its dual space. $B_X$ and $S_X$ denote the closed unit ball and the unit sphere of $X$ respectively. Given two Banach spaces $X$ and $Y$, $\BL(X,Y)$ stands for the Banach space of bounded linear operators from $X$ to $Y$. If $T\in\BL(X,Y)$, $T^*\in\BL(Y^*,X^*)$ denotes the dual operator of $T$. $c_0$ is the Banach space of all real sequences which converge to $0$, endowed with the supremum norm. All Banach spaces involved are assumed to be real.

\section{The space of $C$-fit sequences}
\label{sec:C-fit}
Given a Banach space $X$, recall that $c_0(X)$ is the vector space of all null sequences of $X$, i.e.\ all sequences in $X$ which converge to $0$. This is a Banach space if we endow any $\mathbf{x}={(x_n)}_{n=1}^\infty\in c_0(X)$ with the norm $\norm{\mathbf{x}}{c_0(X)}=\sup\{\norm{x_n}{X}\,:\,n\geq 1\}$. $\ell_1(X)$ is the space of all sequences $\mathbf{x}={(x_n)}_{n=1}^\infty\subset X$ such that the series $\sum_{n=1}^\infty\norm{x_n}{X}$ converges. This is also a Banach space if we endow any $\mathbf{x}={(x_n)}_{n=1}^\infty\in\ell_1(X)$ with the norm $\norm{\mathbf{x}}{\ell_1(X)}=\sum_{n=1}^\infty\norm{x_n}{X}$. We also recall that the dual space of $c_0(X)$ is isometrically isomorphic to $\ell_1(X^*)$. For any $\mathbf{f}={(f_n)}_{n=1}^\infty\in\ell_1(X^*)$ and $\mathbf{x}={(x_n)}_{n=1}^\infty\in c_0(X)$, the evaluation of $\mathbf{f}$ at $\mathbf{x}$ is given by $\mathbf{f}(\mathbf{x})=\sum_{n=1}^\infty f_n(x_n)$. Finally, we call $c_{00}(X)$ the subspace of $c_0(X)$ and $\ell_1(X)$ consisting of all sequences ${(x_n)}_{n=1}^\infty$ such that $x_n\ne 0$ only for finitely many indices $n$. $c_{00}(X)$ is dense in both $c_0(X)$ and $\ell_1(X)$.

The idea for what follows next comes from a well-known characterisation of compactness in Banach spaces due to Grothendiek, which is commonly referred to as the Grothendiek compactness principle: a closed subset $K$ of a Banach space $X$ is compact if and only there is a sequence $\mathbf{x}\in c_0(X)$ such that $K\subseteq\overline{\textup{conv}}(\mathbf{x})$. The Grothendiek compactness principle will turn out to be useful in combination with the main result of \cite{daverava1} (Theorem 4.3), which provides a characterisation of compactivorous sets in Banach spaces. We are especially interested in the second assertion of the theorem.
\begin{theo}
\label{theo:daverava}
Let $E$ be a subset in a Banach space $X$. Then the following assertions are equivalent. 
\begin{enumerate}
\item $E$ is compactivorous.
\item For every compact set $K\subset X$, there are $x\in X$ and $\delta>0$ such that $x+\delta K\subset E$.
\item There is $r>0$ with the property that, for every compact set $K\subset rB_X$, one can find an $x\in X$ such that $x+K\subset E$.
\end{enumerate}
\end{theo}

Since we want to translate compact sets into closed, convex sets, the Grothendiek compactness principle tells us that it is enough to focus on null sequences: given a closed, convex set $C$ in a Banach space $X$, $C$ is compactivorous if and only if for any $\mathbf{x}\in c_0(X)$ there are $\delta>0$ and $z\in X$ such that $z+\delta\mathbf{x}\subset C$.

We say that a sequence $\mathbf{x}={(x_n)}_{n=1}^\infty\in c_0(X)$ is \emph{symmetric} if $x_{2n-1}=-x_{2n}$ for every $n$. Define the \emph{symmetrisation operator} $\func{\Sigma}{c_0(X)}{c_0(X)}$ in the following way: given $\mathbf{x}={(x_n)}_{n=1}^\infty\in c_0(X)$, set $\Sigma(\mathbf{x})=\mathbf{y}$, where $\mathbf{y}={(y_n)}_{n=1}^\infty$ is given by
\[ y_n=\left\{\begin{array}{l l}
		x_k & \text{ if }n=2k-1 \\
		-x_k & \text{ if }n=2k
		\end{array}\right.. \]
In simpler words, a sequence $(x_1,x_2,\dots)$ is mapped to $(x_1,-x_1,x_2,-x_2,\dots)$. Clearly, $\Sigma$ is a linear isometry and $\Sigma(\mathbf{x})$ is symmetric for every $\mathbf{x}\in c_0(X)$. Since $\mathbf{x}$ is a subsequence of $\Sigma(\mathbf{x})$ for every $\mathbf{x}\in c_0(X)$, we can further refine our previous observation: a closed, convex set $C$ in a Banach space $X$ is compactivorous if and only if for any $\mathbf{x}\in c_0(X)$ there are $\delta>0$ and $z\in X$ such that $z+\delta\Sigma(\mathbf{x})\subset C$. In what follows, it will be much more convenient to work with rescalings of $C$ rather than rescalings of sequences.

Let $C$ be a closed, bounded and convex set in a Banach space $X$. We say that a sequence $\mathbf{x}\in c_0(X)$ is \emph{$C$-fit} if there are $z\in X$ and $\lambda>0$ such that $z+\Sigma(\mathbf{x})\subset \lambda C$. If $\mathbf{x}$ if a $C$-fit sequence, its \emph{$C$-fitting cost} is the value
\[ \phi_C(\mathbf{x})=\inf\bigl\{\norm{z}{X}+\lambda\,:\,z\in X,\,\lambda>0,\,z+\Sigma(\mathbf{x})\subset\lambda C\bigr\}. \]
We denote by $F_C(X)$ the subset of $c_0(X)$ consisting of all $C$-fit sequences. Clearly, $F_C(X)=c_0(X)$ if and only if $C$ is compactivorous.

\begin{theo}
Let $X$ be a Banach space and $C\subset X$ be a closed, convex and bounded set. Then $F_C(X)$ is a vector subspace of $c_0(X)$ and $\phi_C$ is a norm on $F_C(X)$ which turns $F_C(X)$ into a Banach space. Moreover, the inclusion operator $\func{I}{F_C(X)}{c_0(X)}$ is continuous.
\end{theo}
\begin{proof}
Let $\mathbf{x}$ and $\mathbf{y}$ be two $C$-fit sequences and let $\epsilon>0$. Let $z_1,z_2\in X$, $\lambda_1>0$ and $\lambda_2>0$ be such that
\begin{align*}
&\norm{z_1}{X}+\lambda_1<\phi_C(\mathbf{x})+\epsilon/2,\quad z_1+\Sigma(\mathbf{x})\subseteq\lambda_1C, \\
&\norm{z_2}{X}+\lambda_2<\phi_C(\mathbf{y})+\epsilon/2,\quad z_2+\Sigma(\mathbf{y})\subseteq\lambda_2C.
\end{align*}  
Then it holds that
\begin{gather*} 
\norm{z_1+z_2}{X}+\lambda_1+\lambda_2<\phi_C(\mathbf{x})+\phi_C(\mathbf{y})+\epsilon, \\
z_1+z_2+\Sigma(\mathbf{x}+\mathbf{y})\subseteq(\lambda_1+\lambda_2)C.
\end{gather*}
This shows that $\mathbf{x}+\mathbf{y}$ is $C$-fit and that $\phi_C(\mathbf{x}+\mathbf{y})<\phi_C(\mathbf{x})+\phi_C(\mathbf{y})+\epsilon$. Since $\epsilon$ is arbitrary, we have the inequality $\phi_C(\mathbf{x}+\mathbf{y})\leq\phi_C(\mathbf{x})+\phi_C(\mathbf{y})$ for all $\mathbf{x},\mathbf{y}\in F_C(X)$.

Let $\mathbf{x}$ be a $C$-fit sequence and let $t>0$. Choose $\epsilon>0$. Find $z\in X$ and $\lambda>0$ such that $\norm{z}{X}+\lambda<\phi_C(\mathbf{x})+\epsilon/t$ and $z+\Sigma(\mathbf{x})\subseteq\lambda C$. Then $\norm{tz}{X}+t\lambda<t\phi_C(\mathbf{x})+\epsilon$ and $tz+\Sigma(t\mathbf{x})\subseteq t\lambda C$. This shows that $t\mathbf{x}$ is $C$-fit and that $\phi(t\mathbf{x})<t\phi_C(\mathbf{x})+\epsilon$. Since $\epsilon$ is arbitrary, we have the inequality $\phi(t\mathbf{x})\leq t\phi_C(\mathbf{x})$ for all $\mathbf{x}\in F_C(X)$ and $t>0$. To see the opposite inequality, just notice that $\phi_C(\mathbf{x})=\phi_C(t^{-1}t\mathbf{x})\leq t^{-1}\phi_C(t\mathbf{x})$. If $t<0$, observe that, since the sequence $\Sigma(\mathbf{x})$ is symmetric, given $z\in X$ and $\lambda>0$, we have $z+\Sigma(\mathbf{x})\subseteq\lambda C$ if and only if $z+\Sigma(-\mathbf{x})\subseteq\lambda C$. This implies that $\phi_C(\mathbf{x})=\phi_C(-\mathbf{x})$. Hence, $\phi_C(t\mathbf{x})=\phi_C(-\lvert t\rvert\mathbf{x})=\lvert t\rvert\phi_C(-\mathbf{x})=\lvert t\rvert\phi_C(\mathbf{x})$, as wished. 

Now, define $R=\sup\{\norm{x}{X}\,:\,x\in C\}$ and let $\mathbf{x}={(x_n)}_{n=1}^\infty$ be a $C$-fit sequence. Choose $\epsilon>0$. Find $z\in X$ and $\lambda>0$ such that $\norm{z}{X}+\lambda<\phi_C(\mathbf{x})+\epsilon$ and $z+\Sigma(\mathbf{x})\subseteq\lambda C$. For each $n$, we have $z+x_n\in\lambda C$ and by the triangular inequality of $X$ we get 
\[\norm{x_n}{X}\leq\norm{z}{X}+\lambda R\leq\max\{1,R\}\cdot(\norm{z}{X}+\lambda)<\max\{1,R\}\cdot\bigl(\phi_C(\mathbf{x})+\epsilon\bigr). \]
Since $\epsilon$ is arbitrary, we have $\norm{x_n}{X}\leq\max\{1,R\}\cdot\phi_C(\mathbf{x})$ for every $n$, which means that $\norm{\mathbf{x}}{c_0(X)}\leq\max\{1,R\}\cdot\phi_C(\mathbf{x})$. This concludes the proof that $\phi_C$ is a norm: if $\phi_C(\mathbf{x})=0$, then $\norm{\mathbf{x}}{c_0(X)}=0$, therefore $\mathbf{x}=\mathbf{0}$. Moreover, we have just shown that the inclusion operator $\func{I}{F_C(X)}{c_0(X)}$ is continuous. 

It remains to check that $\phi_C$ is complete. Let ${(\mathbf{x}_n)}_{n=1}^\infty$ be a Cauchy sequence in $F_C(X)$. By the continuity of $I$, ${(\mathbf{x}_n)}_{n=1}^\infty$ is a Cauchy sequence in $c_0(X)$ as well and it converges to some $\mathbf{x}={(x_j)}_{j=1}^\infty\in c_0(X)$. Let ${(\mathbf{y}_k)}_{k=1}^\infty$ be a subsequence of ${(\mathbf{x}_n)}_{n=1}^\infty$ such that $\phi_C(\mathbf{y}_{k+1}-\mathbf{y}_k)<2^{-k}$ for every $k$ and write $\mathbf{y}_k={(y_{k,j})}_{j=1}^\infty$. Choose $\epsilon>0$ and let $m$ be a positive integer such that $\sum_{k=m}^\infty 2^{-k}<\epsilon$. For each $k$, find $z_k\in X$ and $\lambda_k>0$ such that $\norm{z_k}{X}+\lambda_k<2^{-k}$ and $z_k+\Sigma(\mathbf{y}_{k+1}-\mathbf{y}_k)\subseteq\lambda_kC$. Define \[ z=\sum_{k=m}^\infty z_k,\quad \lambda=\sum_{k=m}^\infty \lambda_k. \] 
For every $k$ and $j$ there is $w_{k,j}\in C$ such that $z_k+y_{k+1,j}-y_{k,j}=\lambda_k w_{k,j}$. Since $y_{k,j}\to x_j$ in $X$ as $k\to\infty$, we get
\[ z+x_j-y_{m,j}=\sum_{k=m}^\infty(z_k+y_{k+1,j}-y_{k,j})=\lambda\sum_{k=m}^\infty\frac{\lambda_k}{\lambda}w_{k,j}. \]
The last series is an infinite convex combination of elements of $C$ and it represents therefore an element $w_j\in C$, as $C$ is closed, convex and bounded. Hence, we have $z+x_j-y_{m,j}=\lambda w_j\in\lambda C$ for every $j$. Similarly, one can prove that, for every $j$, $z-(x_j-y_{m,j})\in\lambda C$. Thus, $z+\Sigma(\mathbf{x}-\mathbf{y}_m)\subset\lambda C$. This shows that $\mathbf{x}-\mathbf{y}_m$ is $C$-fit, which implies that also $\mathbf{x}=(\mathbf{x}-\mathbf{y}_m)+\mathbf{y}_m$ is $C$-fit. Moreover,
\[ \phi_C(\mathbf{x}-\mathbf{y}_m)\leq\norm{z}{X}+\lambda\leq\sum_{k=m}^\infty\bigl(\norm{z_k}{X}+\lambda_k\bigr)<\sum_{k=m}^\infty 2^{-k}<\epsilon. \]
Since $\epsilon$ is arbitrary, the sequence ${(\mathbf{y}_k)}_{k=1}^\infty$ converges to $\mathbf{x}$ in $F_C(X)$ and so does ${(\mathbf{x}_n)}_{n=1}^\infty$, because ${(\mathbf{x}_n)}_{n=1}^\infty$ is a Cauchy sequence.
\end{proof}

From now on, given a closed, convex and bounded set $C$ in a Banach space $X$, we denote by $\phi_C^*$ the norm of the dual space ${F_C(X)}^*$.

\section{Haar positive closed, convex sets}
\label{sec:hpccs}

We recall that, whenever $E$ is a Haar positive set in a separable Banach space, $E$ enjoys the Steinhaus property: $0\in\textup{int}(E-E)$ (see \cite{benlin}, Proposition 6.4). Another crucial tool we will need is the so-called $(\delta,r)$ condition: a characterisation of Haar positive sets due to E.\ Matou\v{s}kov\'{a} (see \cite{eva2}, Theorem 1.1).
\begin{theo}
\label{theo:eva}
Let $E$ be a Borel set in a separable Banach space $X$. $E$ is Haar positive if and only if there are $\delta>0$ and $r>0$ such that, whenever $\mu$ is a Borel probability measure with support in $rB_X$, it is possible to find $x\in X$ such that $\mu(x+E)>\delta$.
\end{theo}
We are ready to work towards a proof of the main theorem. We start with a lemma.

\begin{lemma}
\label{lemma:hpccs 2}
Let $X$ be a Banach space. If $C$ is a subset of $X$ with the Steinhaus property, then $C$ norms $X^*$ in the following way: there is a constant $k>0$ such that for every $f\in X^*$ one has
\[ \sup_{x\in C}\lvert f(x)\rvert\geq k\norm{f}{X^*}. \]
\end{lemma}
\begin{proof}
By definition, there is $r>0$ such that $rB_X\subseteq C-C$. Given $f\in X^*$, we have
\[ \sup_{x\in C}|f(x)|\geq\frac{1}{2}\sup_{x\in C-C}|f(x)|\geq\frac{1}{2}\sup_{x\in rB_X}|f(x)|=\frac{r}{2}\norm{f}{X^*}. \]
Thus, it suffices to set $k=r/2$. 
\end{proof}

\begin{theo}
\label{theo: main1}
Let $C$ be a closed, convex and bounded subset of a separable Banach space $X$. If $C$ is Haar positive, then $F_C(X)=c_0(X)$. In particular, $C$ is compactivorous.
\end{theo}
\begin{proof}
There is no loss of generality if we assume that $0\in C$. The goal is to show that the inclusion $\func{I}{F_C(X)}{c_0(X)}$ is onto. Equivalently, we can establish instead that the dual map $\func{I^*}{\ell_1(X^*)}{{F_C(X)}^*}$ is an embedding, i.e.\ that there is $\gamma>0$ such that $\phi_C^*\bigl(I^*(\mathbf{f})\bigr)>\gamma\norm{\mathbf{f}}{\ell_1(X^*)}$ for all $\mathbf{f}\in\ell_1(X^*)$. By normalisation and density, it suffices to find $\gamma$ such that $\phi_C^*\bigl(I^*(\mathbf{f})\bigr)>\gamma$ for all $\mathbf{f}\in c_{00}(X^*)\cap S_{\ell_1(X^*)}$. 

We start by defining some constants. First, we set $R=\sup\{\norm{x}{X}\,:\,x\in C\}$. Since $-C$ is Haar positive, there are $\delta>0$ and $r>0$ such that $-C$ fulfills the $(\delta,r)$ condition of Theorem \ref{theo:eva}: for every Borel probability measure $\mu$ with support in $rB_X$, there is $z\in X$ such that $\mu(z-C)>\delta$. We may assume that $r\leq R$, so that $rR^{-1}\leq 1$. Finally, by Lemma \ref{lemma:hpccs 2} there is $k>0$ such that, for all $f\in X^*$, 
\[ \sup_{x\in C}\lvert f(x)\rvert\geq k\norm{f}{X^*}. \]

Let $\mathbf{f}={(f_j)}_{j=1}^\infty\in c_{00}(X^*)$ be such that $\norm{\mathbf{f}}{\ell_1(X^*)}=1$ and let $n$ be the minimal index such that $f_j=0$ for all $j\geq n+1$. For each $j\in\{1,\dots,n\}$, find $x_j\in C$ such that $\lvert f_j(x_j)\rvert\geq(k/2)\cdot\norm{f_j}{X^*}$ and let $\mu_j$ be the Dirac probability measure concentrated on $y_j=rR^{-1}x_j\in C$. Set
\[ \mu=\sum_{j=1}^n\norm{f_j}{X^*}\mu_j. \] 
Since $\mu$ is a Borel probability measure whose support is contained in $rB_X$, there must be $z\in X$ such that $\mu(z-C)>\delta$. In other terms,
\[ \sum_{z-y_j\in C}\norm{f_j}{X^*}>\delta. \]
We define $\mathbf{y}'={(y_j')}_{j=1}^\infty\in c_{00}(X)$ in the following way:
\[ y_j'=\left\{\begin{array}{l l}
		\textup{sgn}f_j(y_j)\cdot y_j & \text{ if }z-y_j\in C \\
		0 & \text{ if }z-y_j\notin C\text{ or }j\geq n+1
		\end{array}
	\right.. \]
We now want to check that $\mathbf{y}'$ is $C$-fit. First, notice that $z\in C+y_{j_0}$ for some $j_0\in\{1,\dots,n\}$, hence $z\in 2C$. Every term of the sequence $\Sigma(\mathbf{y}')$ is either $0$, $-y_j$ or $+y_j$ for some $j$ such that $z-y_j\in C$. In the first case we have $z+0=z\in 2C\subseteq 3C$, in the second case we have $z-y_j\in C\subseteq 3C$ and finally, in the third case, we have
\[ z+y_j=z-y_j+2y_j\in C+2C=3C. \]
This shows that $z+\Sigma(\mathbf{y}')\subset 3C$, hence $\mathbf{y}'\in F_C(X)$, as wished. Moreover, we have the estimate $\phi_C(\mathbf{y}')\leq\norm{z}{X}+3\leq 2R+3$. 

Now we have
\[ \delta < \sum_{z-y_j\in C}\norm{f_j}{X^*}\leq\sum_{z-y_j\in C}\frac{2\lvert f_j(x_j)\rvert}{k}=\sum_{z-y_j\in C}\frac{2R\lvert f_j(y_j)\rvert}{rk}, \]
which implies
\begin{align*} 
\frac{\delta rk}{2R} & < \sum_{z-y_j\in C}\lvert f_j(y_j)\rvert=\sum_{j=1}^nf_j(y_j')=\mathbf{f}\bigl(I(\mathbf{y}')\bigr)= \\
                     &=I^*(\mathbf{f})(\mathbf{y}')\leq\phi_C^*\bigl(I^*(\mathbf{f})\bigr)\cdot\phi_C(\mathbf{y}')\leq(2R+3)\cdot\phi_C^*\bigl(I^*(\mathbf{f})\bigr). 
\end{align*}
Therefore, it is enough to set
\[ \gamma=\frac{\delta rk}{2R(2R+3)} \] 
to obtain $\phi_C^*\bigl(I^*(\mathbf{f})\bigr)>\gamma$ for all $\mathbf{f}\in c_{00}(X^*)\cap S_{\ell_1(X^*)}$.
\end{proof}

Theorem \ref{theo: main1} extends in a straightforward way to unbounded sets and delivers the full result we were aiming to.
\begin{theo}
\label{theo: main2}
Let $C$ be a closed, convex set in a separable Banach space $X$. Then $C$ is Haar positive if and only if $C$ is not Haar meagre. Equivalently, $C$ is Haar null if and only if $C$ is Haar meagre.
\end{theo}
\begin{proof}
We have already observed after Proposition \ref{prop:intro} that if $C$ is not Haar meagre, then it must be Haar positive. Conversely, assume that $C$ is Haar positive and write
\[ C=\bigcup_{n=1}^\infty C\cap nB_X. \]
The union on the right side is countable, therefore there must be an $n$ such that $C\cap nB_X$ is Haar positive. By Theorem \ref{theo: main1}, $C\cap nB_X$ is compactivorous and so is $C$. Thus, $C$ is not Haar meagre.
\end{proof}

In \cite{eva3}, E.\ Matou\v{s}kov\'{a} proved that every Haar positive closed, convex set in a separable, reflexive Banach space must have nonempty interior. With Theorem \ref{theo: main2} in our hands, this result can be improved and the proof can be made substantially easier.

\begin{cor}
Let $X$ be a separable dual Banach space and let $C\subseteq X$ be $w^*$-closed and convex. If $C$ is Haar positive, then $\textup{int}\,C\ne\varnothing$ in the norm topology. In particular, in a separable, reflexive Banach space every Haar positive closed, convex set has nonempty interior.
\end{cor}
\begin{proof}
We may assume that $C$ is bounded. Otherwise, find $t>0$ such that $tB_X\cap C$ is Haar positive and work with $tB_X\cap C$ instead, as in the proof of Theorem \ref{theo: main2}. Since $C$ is compactivorous, by Theorem \ref{theo:daverava} there is $r>0$ with the property that every compact subset of $rB_X$ can be translated into $C$. Let ${(F_n)}_{n=1}^\infty$ be an increasing sequence of finite subsets of $rB_X$ such that $\bigcup_{n=1}^\infty F_n$ is $w^*$-dense in $rB_X$. For each $n$, find $z_n$ such that $z_n+F_n\subset C$. Since $C$ is $w^*$-compact, we can assume without loss of generality that the sequence ${(z_n)}_{n=1}^\infty$ has a $w^*$-limit $z\in C$. Now, it is not hard to see that $z+rB_X\subseteq C$. Indeed, take $x\in rB_X$ and find a subsequence ${(F_{n_k})}_{k=1}^\infty$ of ${(F_n)}_{n=1}^\infty$ and a sequence ${(x_k)}_{k=1}^\infty$ with $x_k\in F_{n_k}$ for each $k$ and such that $x_k\overset{w^*}{\rightarrow}x$. Then $z_{n_k}+x_k\overset{w^*}{\rightarrow}z+x$. The assertion for the reflexive case follows easily, as every closed, convex subset of a reflexive Banach space is $w^*$-closed.
\end{proof}

Recall that a cone in a real vector space $X$ is a set $Q$ such that $tx\in Q$ whenever $x\in Q$ and $t$ is a positive scalar. It follows immediately from Theorem \ref{theo:daverava} that a cone $Q$ in a Banach space $X$ is compactivorous if and only if it contains a translation of every compact set. Indeed, let $K\subset X$ be compact and find $x\in X$ and $\delta>0$ such that $x+\delta K\subset Q$. Then $\delta^{-1}x+K\subset\delta^{-1}Q=Q$. If $\mathbf{e}={(e_n)}_{n=1}^\infty$ is a Schauder basis in a Banach space $X$ with coordinate functionals ${(e_n^*)}_{n=1}^\infty$, the positive cone of $\mathbf{e}$ is the set $Q^+(\mathbf{e})=\{x\in X\,:\,e_n^*(x)\geq 0\text{ for every }n\}$. Clearly, $Q^+(\mathbf{e})$ is a closed, convex set. It is proved in \cite{emm_2016} (Corollary 3.2) that, up to equivalence, the standard basis of $c_0$ is the only normalised, unconditional Schauder basis whose positive cone contains a translation of every compact set and therefore is not Haar meagre. If we combine this fact with Theorem \ref{theo: main2}, we can give an affirmative answer to the conjecture raised by the authors of \cite{emm_2016} themselves.
\begin{cor}
Let $\mathbf{e}$ be a normalised, unconditional Schauder basis in a Banach space $X$. If $Q^+(\mathbf{e})$ is Haar positive, then $\mathbf{e}$ is equivalent to the standard basis of $c_0$.
\end{cor}

\section{An application to Banach lattices}
\label{sec:lattices}
In a Banach lattice $X$, the positive cone $X^+=\{x\in X\,:\,x\geq 0\}$ is a closed and convex set. We are interested in characterising those separable Banach lattices whose positive cone is Haar positive. Given a Banach lattice $X$, we shall make use of the following notation:
\[ D_X^+=\biggl\{\bigvee_{j=1}^nx_j\,:\,x_j\in B_X\cap X^+, n\in\N\setminus\{0\}\biggr\},\quad D_X=\{x\in X\,:\,|x|\in D_X^+\}. \]
Clearly, $D_X^+=X^+\cap D_X$ and $B_X\subseteq D_X$.

\begin{lemma}
\label{lemma:AM}
In a Banach lattice $X$ the set $D_X^+$ enjoys the following properties.
\begin{enumerate}
\item For any $x,y\in D_X^+$, also $x\vee y\in D_X^+$.
\item If $x\in D_X^+$ and $0\leq y\leq x$, then $y\in D_X^+$.
\item $D_X^+$ is convex.
\end{enumerate}
\end{lemma}
\begin{proof}
(1) It follows immediately from the definiton of $D_X^+$.

(2) Let $x\in D_X^+$ and consider $y\in X$ such that $0\leq y\leq x$. Find an integer $n\geq 1$ and $x_1,\dots,x_n\in X^+\cap B_X$ such that $x=x_1\vee\cdots\vee x_n$. Notice that
\[ y=y\wedge x=\bigvee_{j=1}^ny\wedge x_j, \]
and that $y\wedge x_j\leq x_j$ for each $j$, which implies $\norm{y\wedge x_j}{X}\leq\norm{x_j}{X}\leq 1$. Thus, $y\wedge x_j\in X^+\cap B_X$ for each $j$, which means that $y\in D_X^+$.

(3) Pick $x,y\in D_X^+$ and $t\in[0,1]$. We have
\[ (1-t)x+ty\leq(1-t)(x\vee y)+t(x\vee y)=x\vee y. \]
By property (1), $x\vee y\in D_X^+$ and, by property (2), $(1-t)x+ty\in D_X^+$. This shows that $D_X^+$ is convex.
\end{proof}

Recall that a Banach lattice $X$ is called an AM-space if, for every $x,y\in X^+$, the equality $\norm{x\vee y}{X}=\max\{\norm{x}{X},\norm{y}{X}\}$ holds. Our last result states that, up to lattice isomorphisms, only separable AM-spaces have nonnegligible positive cones.

\begin{theo}
\label{theo:AM}
Let $X$ be a separable Banach lattice. The following assertions are equivalent.
\begin{enumerate}
\item $X^+$ is Haar positive.
\item For every sequence $\mathbf{x}={(x_n)}_{n=1}^\infty\in c_0(X)$ there exists $z\in X^+$ such that $z\geq\lvert x_n\rvert$ for each $n$.
\item There is $k>0$ (in fact, $k\geq 1$) such that for every integer $n\geq 1$ and any $x_1,\dots,x_n\in X^+$ one has
\[ {\biggl\lvert\biggl\lvert \bigvee_{j=1}^nx_j \biggr\rvert\biggr\rvert}_X\leq k\cdot\max_{1\leq j\leq n}\norm{x_j}{X}. \]
\item There are an AM-space $Y$ and a lattice isomorphism $\func{T}{X}{Y}$.
\end{enumerate}
\end{theo}
\begin{proof}
(1)$\Leftrightarrow$(2) Suppose that $X^+$ is Haar positive and pick $\mathbf{x}={(x_n)}_{n=1}^\infty\in c_0(X)$. By Theorem \ref{theo: main2}, $X^+$ is not Haar meagre and therefore compactivorous, hence there must be $z\in X$ such that $z+\Sigma(\mathbf{x})\subset X^+$. This implies that, for each $n$, we have $z\geq x_n$ and $z\geq -x_n$, hence $z\geq x_n\vee(-x_n)=\lvert x_n\rvert$. Conversely, suppose that (1) holds and pick again $\mathbf{x}={(x_n)}_{n=1}^\infty\in c_0(X)$. Find $z\in X^+$ such that $z\geq\lvert x_n\rvert$ for each $n$. Then $z+\Sigma(\mathbf{x})\subset X^+$. This shows that $X^+$ is compactivorous, therefore Haar positive.

(2)$\Leftrightarrow$(3) Suppose that $X$ fulfils assertion (3) and let $k>0$ be such that 
\[ {\biggl\lvert\biggl\lvert \bigvee_{j=1}^nx_j \biggr\rvert\biggr\rvert}_X\leq k\cdot\max_{1\leq j\leq n}\norm{x_j}{X} \] 
for every integer $n\geq 1$ and any $x_1,\dots,x_n\in X^+$. Fix $\mathbf{x}={(x_n)}_{n=1}^\infty\in c_0(X)$. Define the sequence $\mathbf{y}={(y_n)}_{n=1}^\infty$ recursively in the following way: $y_1=\lvert x_1\rvert$ and $y_{n+1}=y_n\vee\lvert x_{n+1}\rvert$ for every $n$. Observe that $\mathbf{y}$ is an increasing sequence in $X^+$. Choose $\epsilon>0$ and let $n_0$ be such that $\norm{x_n}{X}<\epsilon/k$ for all $n\geq n_0$. For every $n$ and $m$ with $n>m\geq n_0$ we have
\[ 0\leq y_n-y_m=\bigvee_{j=1}^n\lvert x_j\rvert-\bigvee_{j=1}^m\lvert x_j\rvert\leq\bigvee_{j=m+1}^n\lvert x_j\rvert, \]
which implies
\[ \norm{y_n-y_m}{X}\leq{\biggl\lvert\biggl\lvert\bigvee_{j=m+1}^n\lvert x_j\rvert\biggr\rvert\biggr\rvert}_X\leq k\cdot\max_{m+1\leq j\leq n}\norm{x_j}{X}<\epsilon. \]
That is, $\mathbf{y}$ is a Cauchy sequence and it has therefore a limit $z\in X^+$. Clearly, $z\geq y_n$ for every $n$, thus $z\geq\lvert x_n\rvert$ for every $n$. On the other hand, assume by contradiction that for each positive integer $k$ there are $x_{k,1},\dots,x_{k,n(k)}\in B_X\cap X^+$ such that $\norm{x_{k,1}\vee\cdots\vee x_{k,n(k)}}{X}\geq k^2$. Let $\mathbf{y}={(y_n)}_{n=1}^\infty\in c_0(X)$ be the sequence whose elements are $x_{1,1},\dots,x_{1,n(1)}$, followed by $2^{-1}x_{2,1},\dots,2^{-1}x_{2,n(2)}$, followed by $3^{-1}x_{3,1},\dots,3^{-1}x_{3,n(3)}$ and so on. If there were $z\in X^+$ such that $z\geq y_n$ for all $n$, then for each $k$ it would hold that $z\geq k^{-1}x_{k,1}\vee\cdots\vee k^{-1}x_{k,n(k)}$. Thus, for each $k$, 
\[ \norm{z}{X}\geq\norm{k^{-1}x_{k,1}\vee\cdots\vee k^{-1}x_{k,n(k)}}{X}=k^{-1}\norm{x_{k,1}\vee\cdots\vee x_{k,n(k)}}{X}\geq k^{-1}k^2=k, \]
which does not make any sense.

(3)$\Leftrightarrow$(4) Assume that (4) holds, i.e.\ that there are an AM-space $Y$ and a lattice isomorphism $\func{T}{X}{Y}$. Set $k=\norm{T^{-1}}{\BL(Y,X)}\norm{T}{\BL(X,Y)}$ and choose an integer $n\geq 1$ and $x_1,\dots,x_n\in X^+$. Then
\begin{align*}
{\biggl\lvert\biggl\lvert \bigvee_{j=1}^nx_j \biggr\rvert\biggr\rvert}_X &\leq\norm{T^{-1}}{\BL(Y,X)}\cdot{\biggl\lvert\biggl\lvert\,T\biggl(\,\bigvee_{j=1}^nx_j\biggr)\biggr\rvert\biggr\rvert}_Y= \\
									    &=\norm{T^{-1}}{\BL(Y,X)}\cdot{\biggl\lvert\biggl\lvert\bigvee_{j=1}^nT(x_j)\biggr\rvert\biggr\rvert}_Y= \\
									    &=\norm{T^{-1}}{\BL(Y,X)}\cdot\max_{1\leq j\leq n}\norm{T(x_j)}{Y}\leq \\
									    &\leq\norm{T^{-1}}{\BL(Y,X)}\norm{T}{\BL(X,Y)}\cdot\max_{1\leq j\leq n}\norm{x_j}{X}= \\
									    &= k\cdot\max_{1\leq j\leq n}\norm{x_j}{X}.
\end{align*}
To prove the converse implication, we will exhibit an equivalent lattice norm on $X$ which turns $X$ into an AM-space. If (3) holds, then it follows from the definition of $D_X^+$ and from the properties of lattice norms that $\norm{x}{X}\leq k$ for every $x\in D_X$, which implies $B_X\subseteq D_X\subseteq kB_X$. $D_X$ is symmetric, as $|x|=|-x|$ for every $x\in X$. Moreover, $D_X$ is convex. Indeed, if $x,y\in D_X$ and $t\in[0,1]$, then $|x|,|y|\in D_X^+$ and $|(1-t)x+ty|\leq (1-t)|x|+t|y|$. By the convexity of $D_X^+$ and by (2) in Lemma \ref{lemma:AM}, $|(1-t)x+ty|\in D_X^+$, hence $(1-t)x+ty\in D_X$. Let $\nu$ be the Minkowski functional associated to the set $D_X$. That is, $\nu(x)=\inf\{\lambda>0\,:\,x\in\lambda D_X\}$ for every $x\in X$. Since $B_X\subseteq D_X\subseteq kB_X$, $\nu$ defines an equivalent norm on $X$. Moreover, one notices easily that $\nu(x)=\nu(|x|)$ for every $x\in X$. To see that $\nu$ is indeed a lattice norm, pick $x,y\in X$ such that $|x|\leq|y|$. For every $\epsilon>0$, we have
\[ \frac{|x|}{\nu(y)+\epsilon}\leq\frac{|y|}{\nu(y)+\epsilon}=\frac{|y|}{\nu(|y|)+\epsilon}\in X^+\cap D_X=D_X^+. \]
By (2) in Lemma \ref{lemma:AM}, we obtain $|x|\in\bigl(\nu(y)+\epsilon\bigr)D_X^+$. As $\epsilon$ is arbitrary, it follows that $\nu(x)=\nu(|x|)\leq\nu(y)$, as wished. Now, let $x,y\in X^+$ and set $m=\max\{\nu(x),\nu(y)\}$. For every $\epsilon>0$, we have $x,y\in(m+\epsilon)D_X^+$. By (1) in Lemma \ref{lemma:AM}, it follows that
\[ \frac{x\vee y}{m+\epsilon}=\frac{x}{m+\epsilon}\vee\frac{y}{m+\epsilon}\in D_X^+. \]
That is, $x\vee y\in(m+\epsilon)D_X^+$. Since $\epsilon$ is arbitrary, this allows to conclude that $\nu(x\vee y)\leq m$. The opposite inequality is true in every Banach lattice, therefore $\nu$ turns $X$ into an AM-space.
\end{proof}

Remark: the proof of the equivalence between (2), (3) and (4) in Theorem \ref{theo:AM} does not make use of any separability condition. The equivalence holds therefore in general Banach lattices.

\section*{Acknowledgements}
The author would like to thank Professor Eva Kopeck\'{a} for the many motivating discussions and the referee, who helped simplify the arguments a lot. The work of the author is supported by the Austrian Science Fund (FWF): P 32523-N.

\printbibliography
\end{document}